\documentclass[10pt]{amsart}
\usepackage{amssymb}
\usepackage{tikz}
\usepackage{graphicx}
\usepackage{tikz-cd}

\newtheorem{theorem}{Theorem}[section]
\newtheorem{lemma}[theorem]{Lemma}

\theoremstyle{definition}

\theoremstyle{remark}

\numberwithin{equation}{section}

\begin{document}

\begin{center}
\textbf{\Large{Products and connected sums of spheres as monotone Lagrangian submanifolds}}
\end{center}
\begin{center}
\textbf{Vardan Oganesyan and Yuhan Sun}
\end{center}

\begin{abstract}
We obtain new restrictions on Maslov classes of monotone Lagrangian submanifolds of $\mathbb{C}^n$. We also construct families of new examples of monotone Lagrangian submanifolds, which show that the restrictions on Maslov classes are sharp in certain cases.
\end{abstract}

\tableofcontents

\section{Introduction}

\footnote{This work was supported by the Russian Science Foundation under grant no.18-11-00316.}

Let $L$ be a closed Lagrangian submanifold of $\mathbb{C}^{n}$. We can ask the following question:

\begin{center}
\emph{ What can be said about the minimal Maslov number $N_L$ of $L$?  }
\end{center}

\vspace*{.1in}
The famous  Audin conjecture says the following:
\\
\\
\textbf{Conjecture}. Any emnedded Lagrangian torus in $\mathbb{C}^{n}$ has minimal Maslov number $N=2$.
\\

Oh proved the conjecture for monotone tori of dimension less than $24$ (see \cite{Oh}). Buhovsky in \cite{B}, Fukaya-Oh-Ohta-Ono in \cite{FOOO}, Damian in \cite{D} proved the conjecture for monotone tori of arbitrary dimension. The conjecture was proved in full generality by Cieliebak-Mohnke in \cite{CM} and by Irie in \cite{I}. Outside the torus case, there are various restrictions on the Maslov class coming from the Floer theory. One of them is given by Oh's spectral sequence. Oh  proved the following two theorems:
\\

\textbf{Theorem.}(Oh, \cite{Oh})
Let $N_L$ be the minimal Maslov number of a closed monotone Lagrangian submanifold $L \subset \mathbb{C}^{n}$, where $n>1$. Then
\begin{equation*}
1 \leqslant N_L \leqslant n.
\end{equation*}

\textbf{Theorem.}(Oh, \cite{Oh})
Let $n$ be a positive even number and let $N_L$ be the minimal Maslov number of a monotone Lagrangian embedding of $L = S^{n-1}\times S^{1}$ into $\mathbb{C}^{n}$ . Then $N_L$ divides $n$.
\\

In \cite{P} Polterovich constructed monotone Lagrangian embeddings of $S^{n-1}\times S^{1}$ into $\mathbb{C}^n$ with minimal Maslov numbers $2$ and $n$, where $n$ is an even number. However, examples of Lagrangian embeddings $L = S^{n-1}\times S^{1} \rightarrow \mathbb{C}^n$ with minimal Maslov number $N_L  \neq 2, n $ are not known. Also, many interesting results on Lagrangian intersection theory are proved in \cite{BC}.

In papers \cite{MironovCn, MirPan, kotel} Mironov, Panov  and Kotelskiy studied Hamiltonian-minimal and minimal Lagrangian submanifolds of toric manifolds. In particular, they associated a closed Hamiltonian-minimal Lagrangian submanifold  $L \subset \mathbb{C}^n$ to each Delzant polytope (see Sections $\ref{spolytope}$ and $\ref{seclagr}$ for details).  The Lagrangian $L$ is diffeomorphic to the total space of  fiber bundle over $T^k$, where the fiber is the so-called real moment-angle manifold associated to $P$. In some cases real moment angle manifolds are diffeomorphic to connected sums of sphere products, but in general real moment-angle manifolds define a rich family of smooth manifolds and their topology is far from being completely understood (see \cite{Lopez, real1, real2}).
\\
\\
\textbf{Remark.} In fact, Mironov, Panov and Kotelskiy associated a Lagrangian submanifold to a representation of a Delzant polytope by inequalities. It may happen that some of the inequalities can be removed from the presentation without changing the polytope. Such inequalities are called redundant. A presentation without redundant inequalities is called irredundant. It turns out that two systems of inequalities, that provides the same polytope, can give different Lagrangian submanifolds. We will see how it works in Theorem $\ref{exist2}$.
\\

It was noticed by the first author that methods of Mironov, Panov  and Kotelskiy can be used for constructing monotone Lagrangian submanifolds. The following theorem was proved:
\\
\\
\textbf{Theorem.} (\cite{Og}).
\emph{Let $P$ be a Delzant and irredundant polytope and $L \subset \mathbb{C}^n$ be the corresponding Lagrangian submanifold. Then $L$ is monotone if and only if the polytope $P$ is Fano.}
\\

There are many Delzant polytopes (infinitely many even in dimension $2$). Also, there are many Delzant polytopes that are Fano. As a result we obtain a large family of monotone Lagrangian submanifolds. Moreover, the Maslov class of a monotone Lagrangian associated to a Delzant Fano polytope can be found explicitly. It turns out that many polytopes provide the same diffeomorphism type of $L$. As a result we can construct monotone Lagrangian submanifolds of $\mathbb{C}^n$ with different minimal Maslov number, and therefore distinct up to Lagrangian isotopy. Also, it can be shown that some of our embeddings are smoothly isotopic but not Lagrangian isotopic (see \cite{Og} for more details). On the other hand, we have smoothly isotopic Lagrangians with equal minimal Maslov numbers and the Maslov class doesn't allow us to distinguish these Lagrangians up to Lagrangian or Hamiltonian isotopy. By studying many examples of the constructed Lagrangians, we pose the following \textquotedblleft Lagrangian version of Delzant theorem\textquotedblright.
\\
\\
\textbf{Conjecture.} Let $P_1$, $P_2$ be Delzant Fano and irredundant  polytopes in $\mathbb{R}^k$. Let $L_1$, $L_2$ be the corresponding embedded monotone Lagrangian submanifolds  of $\mathbb{C}^n$. Then $L_1$ is Hamiltonian isotopic to $L_2$ if and only if $P_1 = g\cdot P_2$, where $g \in SL(k,\mathbb{Z}) \ltimes translations$.
\\
\\
Removing monotonicity assumption results in a stronger conjecture, which we also expect to be hold:
\\
\\
\textbf{Conjecture.} Let $P_1$, $P_2$ be  Delzant irredundant polytopes in $\mathbb{R}^k$. Let $L_1$, $L_2$ be the corresponding embedded Lagrangian submanifolds of $\mathbb{C}^n$. Then $L_1$ is Hamiltonian isotopic to $L_2$ if and only if $P_1 = g\cdot P_2$, where $g \in SL(k,\mathbb{Z}) \ltimes translations$.
\\

In this paper, we find new restrictions on the Maslov classes of monotone Lagrangian submanifolds of $\mathbb{C}^n$. Then we construct new examples of monotone Lagrangians, which show that the restrictions on Maslov classes are sharp in certain cases. We now proceed to precise formulations of our results.  We prove the following two theorems using lifted Floer homology:

\begin{theorem}\label{rest1}
Let $L=S^{p-1}\times S^{n-p-1}\times T^{2}$ be a monotone Lagrangian submanifold of $\mathbb{C}^{n}$, where $p, n$ are arbitrary positive even integers. Then the minimal Maslov number $N_L$ of $L$ divides either $p$, or $n-p$.
\end{theorem}

\begin{theorem}\label{rest2}
Let $p$ be an arbitrary positive even integer. Assume that $m , l$ are arbitrary positive integers. If $L = (S^{p-1} )^m \times (S^{1})^{l}$ is embedded into $\mathbb{C}^{m(p-1) + l}$ as a monotone Lagrangian submanifold, then the minimal Maslov number $N_L$ divides $p$.

\end{theorem}
\textbf{Remark}. Both authors acknowledge the anonymous  referee for his comments. The referee advised us to use lifted Floer homology to prove these theorems.
\\

One can ask  whether these restrictions on the Maslov class are sharp. To answer this question we construct families of explicit examples. Denote by $gcd(a,b)$ the greatest common divisor of $a$ and $b$.

\vspace{.1in}

\begin{theorem}\label{exist1}
Let $n, \; p, \;k$  be  even positive numbers such that $n-p+k>p$, $k<p-1$, $p>2$, and $n-p>2$. Then there exists a monotone Lagrangian embedding of $S^{p-1} \times S^{n-p-1} \times T^2$ into $\mathbb{C}^n$ with minimal Maslov number $\gcd(p, n-p+k)$.

On the other hand, if $N$ is an arbitrary even divisor of $p$  and $n \geqslant 2p$, then $N$ can be the minimal Maslov number of a Lagrangian monotone embedding of $S^{p-1} \times S^{n-p-1} \times T^2$ into $\mathbb{C}^n$.
\end{theorem}

The Fano polytope associated to the monotone Lagrangian $ S^{p-1}\times S^{n-p-1} \times T^{2}$ is the product of two simplices $\Delta^{p-1}$ and $\Delta^{n-p-1}$. Parameter $k$ depends on the angle between the simplices.

Let us assume that $n=2p$ and $N$ is an arbitrary even divisor of $p$. Then from Theorem $\ref{exist1}$ we get monotone Lagrangian embedding of $S^{p-1} \times S^{p-1} \times T^2$ into $\mathbb{C}^{2p}$ with minimal Maslov number $N$. Also, example of Polterovich gives us a monotone Lagrangian $S^{p-1}\times S^{1} \subset \mathbb{C}^p$ with minimal Maslov number $p$. Taking product of these two families of examples we obtain monotone Lagrangians $(S^{p-1}\times S^{1})^{m}$ with minimal Maslov number $N$, where $m \geqslant 2$. So, we proved the following:
\\
\\
\textbf{Corollary.} Let $p$ be a positive even integer and $N$ be an arbitrary positive even divisor of $p$. Assume that $m > 1$ is an arbitrary integer. Then there exists a monotone Lagrangian embedding of $(S^{p-1}  \times S^{1})^{m}$ into $\mathbb{C}^{mp}$ with minimal Maslov number $N$.
\\

We also want to note that our construction sometimes provides more than one monotone Lagrangian with minimal Maslov number $N$. So, it can be interesting to distinguish them by other invariants (see conjectures above).
\\

We can get more results in a similar way.

\begin{theorem}\label{exist2}
Let $n$ be an odd integer and $k$ be an even integer such that $0 \leqslant k \leqslant n-2$, $k > \frac{n-3}{2}$, $n>3$. There exists a monotone Lagrangian embedding of $L = S^{n-2}\times T^{2}$ into $\mathbb{C}^{n}$ with minimal Maslov number $N_L=\gcd(n-1, 2k+2)$.
\end{theorem}

The Fano polytope associated to the Lagrangian $L = S^{n-2}\times T^{2}$ is a simplex, but the representation of the simplex by inequalities is redundant. (see Sections $\ref{seclagr}$ and $\ref{proofexist2}$ for details) In fact, if $n$ is an even number, then the monotone Lagrangian associated to a standard irredundant Fano representation of a simplex is $S^{n-1} \times S^1$.

By Theorem $\ref{rest2}$ the minimal Maslov number $N_L$ of monotone embedded Lagrangian $L = S^{n-2}\times T^{2} \subset \mathbb{C}^n$ divides $n-1$. Note that $k$ is an even number and $2k+2 \equiv 2 \mod 4$. So, we have the following two cases (see Section $\ref{exist2} for the proof)$:
\begin{enumerate}
\item if $n-1$ is divisible by $4$, then any divisor of $n-1$ of the form $4l+2$ can be the minimal Maslov number of a monotone Lagrangian embedding of $S^{n-2} \times T^2$ into $\mathbb{C}^n$, where $l$ is an non negative integer such that $4l+2 < n-1$.
\item if $n-1 \equiv 2 \mod 4$, then any even divisor of $n-1$ less than $n-1$ can be the minimal Maslov number of a monotone Lagrangian embedding of $S^{n-2} \times T^2$ into $\mathbb{C}^n$.
\end{enumerate}
So, either ``half of all'', or ``all but one'' divisors of $n-1$ can be realized.
\\
\\
\textbf{Examples.} Assume that $n=13$. There exist monotone Lagrangian emebddings of $L = S^{11} \times T^2$ into $\mathbb{C}^{13}$ with minimal Maslov number $N_L = gcd(12, 2k+2)$, where $k$ is even and $5 < k \leqslant 11$. If $k=2$ or $10$, then $N_L = 2$. If $k=8$, then $N_L = 6$.

Assume that $n=31$. There exist monotone Lagrangian emebddings of $L = S^{29} \times T^2$ into $\mathbb{C}^{31}$ with minimal Maslov number $N_L = gcd(30, 2k+2 )$, where $k$ is even and $14<k<29$. If $k=16, 18, 22, 28$, then $N_L = 2$. If $k=20$ or $k=26$, then $N_L = 6$. If $k=24$, then $N_L = 10$.

It would be interesting to distinguish embeddings with equal minimal Maslov numbers by other invariants.

We can get similar theorems for Lagrangians with more complicated topology.

\begin{theorem}\label{exist3}
Let $p$ be an even positive integer and $N$ be any positive even divisor of $p$. Then there exists a monotone Lagrangian embedding of $(\#_5 (S^{2p-1} \times S^{3p-2}))\times T^3$ into $\mathbb{C}^{5p}$ with minimal Maslov number $N$.

Moreover, the minimal Maslov number of any monotone Lagrangian embedding of $(\#_5 (S^{2p-1} \times S^{3p-2}))\times T^3$ divides $p$.
\end{theorem}

\subsection*{Acknowledgements}
Both authors acknowledge Kenji Fukaya and Mark McLean for many helpful discussions during the preparation of this note. Also, authors acknowledge the anonymous  referee for useful comments.

\section{Preliminaries}\label{prelim}
Let $\mathbb{C}^{n}$ be the complex vector space with the standard symplectic form $\omega$ and let $L$ be a closed Lagrangian submanifold of $\mathbb{C}^n$. There are two homomorphisms on $\pi_2(\mathbb{C}^n, L)$: the symplectic energy
\begin{equation*}
I_{\omega}: \pi_2(\mathbb{C}^n, L) \rightarrow \mathbb{R}
\end{equation*}
and the Maslov index
\begin{equation*}
I_{\mu}: \pi_2(\mathbb{C}^n, L) \rightarrow \mathbb{Z}.
\end{equation*}
We have that $\pi_2(\mathbb{C}^n, L) = \pi_1(L)$. Hence, we can assume that both homomorphisms are defined on $H_1(L, \mathbb{Z})$, i.e.
\begin{equation*}
\begin{gathered}
I_{\omega}: H_1(L, \mathbb{Z}) \rightarrow \mathbb{R} \\
I_{\mu}: H_1(L, \mathbb{Z}) \rightarrow \mathbb{Z}.
\end{gathered}
\end{equation*}

\textbf{Definition.} We say that $L$ is \textit{monotone} if $I_{\omega}=c I_{\mu}$ for some positive real constant $c$. The \textit{minimal Maslov number} $N_L$ of $L$ is the positive generator of $I_{\mu}(H_1(L, \mathbb{Z})) \subset \mathbb{Z}$.
\\

If we assume that $L$ is orientable, then $I_{\mu}(H_1(L, \mathbb{Z}))$  is contained in $2\mathbb{Z}$ and the minimal Maslov number $N_L$ is even.

\subsection{Lifted Floer Homology and related spectral sequence.} Assume that $L$ is a monotone Lagrangian submanifold of $\mathbb{C}^n$ with $N_L  \geqslant 3$. Let $\tilde{L}$ be the universal cover of $L$. Damian in \cite{D} defined the lifted Floer homology $FH^{\tilde{L}}(L)$. We do not give definition of the lifted Floer homology in this paper, but we want to mention that $FH^{\tilde{L}}(L) = 0$ because $L \subset \mathbb{C}^n$. Moreover, Damian constructed a spectral sequence which converges to the lifted Floer homology and whose first page is built using the singular homology of $\tilde{L}$. Let us give more details. Denote by $A$ a ring $\mathbb{Z}_2[T, T^{-1}]$, where we assume that $deg(T) = N_L$. Also, denote by $A^k \subset A$ a subgroup $\mathbb{Z}_2 \cdot T^k$.
\\
\\
\textbf{Theorem} (Damian \cite{D}). There exists a spectral sequence $(E^{t,s}_r, d_r)$ with the following properties:
\\
\\
$\bullet$ $E^{t,s}_r$ converges to the lifted Floer homology $FH^{\tilde{L}}(L)$. In our case $L \subset \mathbb{C}^n$ and $FH^{\tilde{L}}(L) = 0$.
\\
\\
$\bullet$ $E_1^{t,s} = H_{t+s - tN_L}(\tilde{L}, \mathbb{Z}_2) \otimes A^{t}$, $d_1 = \delta_1 \otimes T^{-1}$, where
\begin{equation*}
\delta_1: H_{t+s - tN_L}(\tilde{L}, \mathbb{Z}_2) \rightarrow H_{t+s - 1- (t-1)N_L }(\tilde{L}, \mathbb{Z}_2).
\end{equation*}
Let us recall that $deg(T) = N_L$.
\\
\\
$\bullet$ For every $r >1$, element $E^{t,s}_{r}$ has the form $E^{t,s}_{r} = V^{t,s}_r \otimes A^t$ with $d_r = \delta_r \otimes T^{-r}$, where  $V^{t,s}_r$ is vector space over $\mathbb{Z}_2$ and $\delta_r : V^{t,s}_r \rightarrow V^{t-r,s+r-1}_r$ is a homomorphism defined for every $t, s$  such that $\delta_r \circ \delta_r = 0$. We have
\begin{equation*}
V^{t,s}_{r+1} = \frac{Ker(\delta_r: V^{t,s}_r \rightarrow V^{t-r,s+r-1}_r)}{Im(\delta_r: V^{t+r,s-r+1}_r \rightarrow V^{t,s}_r)}.
\end{equation*}

$\bullet$ The spectral sequence collapses at page $[\frac{dim(L) + 1 }{N_L} ] +1$.

\section{Polytopes and intersection of quadrics}\label{spolytope}
In this section we discuss toric topology and its applications. For more details we refer our reader to paper of Panov \cite{Panov} (Sections $2, 3, 12$). Much more details can be found in book \cite{torictop}.

A convex \emph{polyhedron} $P$ is an intersection of finitely many halfspaces in $\mathbb{R}^k$. Bounded polyhedra are called \emph{polytopes}.

A \emph{supporting hyperplane} of $P$ is a hyperplane $H$ which has common points with $P$ and for which the polyhedron is contained in one of the two closed half-spaces determined by $H$. The intersection $P\cap H$ with a supporting hyperplane is called a face of the polyhedron. Zero-dimensional faces are called \emph{vertices} and faces of codimension one are called \emph{facets}.

Consider a system of $n$ linear inequalities defining a convex polyhedron in $\mathbb{R}^k$
\begin{equation}\label{polytope}
P_{A,b}=\{x\in \mathbb{R}^k: <a_i,x>+b_i \geqslant 0 \quad for \quad i=1,...,n  \},
\end{equation}
where $<\cdot,\cdot>$ is the standard scalar product on $\mathbb{R}^k$, $a_i \in \mathbb{R}^k$, and $b_i \in \mathbb{R}$. By $b$ denote a vector $b=(b_1,...,b_n)^T$, $x=(x_1,...,x_k)^T$ and by $A$ the $k\times n$ matrix whose columns are the vectors $a_i$. Then, our polyhedron can be written in the following form:
\begin{equation*}
P_{A,b}=\{x\in \mathbb{R}^k: (A^Tx + b)_i \geqslant 0 \quad for \quad i=1,...,n\}.
\end{equation*}

\vspace*{.1in}

\textbf{Definition.}  We say that (\ref{polytope}) is \emph{simple} if exactly $k$ facets meet at each vertex. We say that (\ref{polytope}) is \emph{generic} if for any vertex $x\in P$ the normal vectors $a_i$ of the hyperplanes containing $x$ are linearly independent.
\\

Let us assume that $\mathbb{Z}<a_1,....,a_n>$ defines a lattice in $\mathbb{R}^{n-k}$, where $\mathbb{Z}<a_1,....,a_n>$ is the set of integer linear combinations of vectors $a_1, ..., a_n$.
\\
\\
\textbf{Definition.} Polyhedron $P$ is called \emph{Delzant} if it is simple and for any vertex $x \in P$ the vectors $a_i$ normal to the facets meeting at $x$ form a basis for the lattice $\mathbb{Z}<a_1,....,a_n>$.
\\
\\
\textbf{Definition.} A Delzant polytope is called \emph{Fano} if it can be defined by
\begin{equation*}
P_{A,b}=\{x\in \mathbb{R}^k: \langle a_i,x \rangle + C \geq  0 \quad for \quad i=1,...,n  \},
\end{equation*}
where  each vector $a_i \in \mathbb{Z}^k$ is the primitive integral interior normal to the corresponding facet. In other words, $C=b_1=...=b_n$.
\\

System of inequalities $(\ref{polytope})$ gives us a linear map from $\mathbb{R}^k$ to $\mathbb{R}^n$
\begin{equation}\label{polytope2}
\begin{gathered}
i_{A,b} : \mathbb{R}^k \rightarrow \mathbb{R}^n, \\
i_{A,b}(x)=A^Tx + b = (\langle a_1, x \rangle + b_1,...,\langle a_n, x \rangle + b_n)^T.
\end{gathered}
\end{equation}
Then, the image $i_{A,b}(\mathbb{R}^k)$ can be written as
\begin{equation}\label{polsys}
\begin{gathered}
i_{A,b}(\mathbb{R}^k)= \{u \in \mathbb{R}^n : \Gamma u = \Gamma b \},\\
\Gamma A^T = 0, \quad u=(u_1,...,u_n)^T,
\end{gathered}
\end{equation}
where $\Gamma $ is $(n-k)\times n$-matrix whose rows form a basis of linear relations between the vectors $a_i$. We have
\begin{equation}\label{polytopeemb}
i_{A,b}(P) = i_{A,b}(\mathbb{R}^k)\cap \mathbb{R}^n_{+}.
\end{equation}

Let us describe the correspondence between the intersection of quadrics and polyhedra. Replacing $u_i$ by $u_i^2$ in (\ref{polsys}) we get $(n-k)$ quadrics which define a subset in $\mathbb{R}^n$.

Now assume that we have
\begin{equation}\label{quadrics}
\mathcal{R}_{\Gamma, \delta} = \{u\in \mathbb{R}^n : \gamma_{1,i}u_1^2 + ... + \gamma_{n,i}u_n^2 = \delta_i, \quad i=1,...,n-k, \}.
\end{equation}
The coefficients of the quadrics define $(n-k)\times n$ matrix $\Gamma=(\gamma_{jk})$. The group $\mathbb{Z}_{2}^n$ acts on $\mathcal{R}_{\Gamma,\delta}$ by
\begin{equation*}
\varepsilon \cdot (u_1,...,u_n) = (\pm u_1,...,\pm u_n)
\end{equation*}
The quotient $\mathcal{R}_{\Gamma,\delta}/\mathbb{Z}_{2}^n$ can be identified with the set of nonnegative solutions of the system
\begin{equation*}
\left\{
 \begin{array}{l}
\gamma_{1,1}u_1 + ... + \gamma_{n,1}u_n = \delta_1\\
\ldots \quad \quad \\
\gamma_{1,(n-k)}u_1 + ... + \gamma_{n,(n-k)}u_n = \delta_{n-k}
 \end{array}
\right.
\end{equation*}
Easy to see that the system above is equivalent to (\ref{polsys}) and (\ref{polytopeemb}). Solving the homogeneous version of the system above we get the matrix $A$. We have that rows of matrix $\Gamma$ form a basis of linear relations between the vectors $a_i$. Then, we can construct a polytope (\ref{polytope}), where $b=(b_1,...,b_n)$ is an arbitrary solution of the linear system above.

We obtain that a polyhedron defines an intersection of quadrics and an intersection of quadrics defines a polyhedron.

It may happen that some of the inequalities can be removed from the presentation without changing $P_{A,b}$. Such inequalities are called \emph{redundant}. A presentation without redundant inequalities is called \emph{irredundant}.
\\
\\
\begin{theorem}\label{redundant}(\cite{Panov} Theorem 3.5 and Chapter 12, \cite{MirPan}).

The intersection of quadrics $\mathcal{R}_{\Gamma, \delta}$ is nonempty and nondegenerate if and only if the presentation $P_{A,b}$ is generic.

Assume that we have a polytope $P$ given by two different systems of inequalities. Assume that our inequalities are defined by matrices  $A$ and $A'$. Let $\mathcal{R}_{\Gamma,\delta}$, $\mathcal{R}_{\Gamma',\delta'}$ be the corresponding intersection of quadrics.  If system $A'$ is obtained from $A$ by adding $m$ redundant inequalities, then $\mathcal{R}_{\Gamma',\delta'}$ is homeomorphic to a product of $\mathcal{R}_{\Gamma,\delta}$ and $\mathbb{Z}_2^m$, i.e. $\mathcal{R}_{\Gamma',\delta'}$ is disjoint union of $2^m$ copies of $\mathcal{R}_{\Gamma,\delta}$.

\end{theorem}

\section{Lagrangian submanifolds of $\mathbb{C}^n$}\label{seclagr}
A construction explained in this section was discovered by Mironov in \cite{MironovCn}. Later, the construction was studied by Panov and Kotelskiy in \cite{MirPan} and \cite{kotel}. They applied Mironov's idea to study minimal and Hamiltonian-minimal Lagrangians of toric manifolds. It was noticed by the first author in \cite{Og} that their method can be used to construct monotone Lagrangian submanifolds. Moreover, the construction allows us to find the Maslov class of the Lagrangians.

Let $\mathcal{R}$ be a $k$-dimensional submanifold of $\mathbb{R}^n$ defined by
\begin{equation}\label{eqmain}
\begin{gathered}
\left\{
 \begin{array}{l}
\gamma_{1,i}u_1^2 + ... + \gamma_{n,i}u_n^2 = \delta_i, \quad i=1,...,n-k,
 \end{array}
\right.
\\
\delta_i \in \mathbb{R}, \quad \gamma_{i,j} \in \mathbb{Z}
\end{gathered}
\end{equation}
Let $\gamma_j$ be the $j$th column of system $(\ref{eqmain})$, i.e.
\begin{equation}\label{vectors}
\gamma_j=(\gamma_{j,1},...,\gamma_{j, n-k})^T \in \mathbb{Z}^{n-k}, \quad j=1,...,n
\end{equation}
and $\Gamma$ be a matrix with columns $\gamma_j$, $j=1,...,n$. Let us assume that the integer vectors $\gamma_1$,...,$\gamma_n$ form a lattice $\Lambda$ in $\mathbb{R}^{n-k}$ of maximum rank. The dual lattice $\Lambda^{*}$ is defined by
\begin{equation*}
\Lambda^{*}=\{\lambda^{*} \in \mathbb{R}^{n-k}| \langle\lambda^{*},\lambda\rangle \in \mathbb{Z}, \; \lambda \in \Lambda\},
\end{equation*}
where $\langle \lambda^{*},\lambda \rangle$ is the Euclidian product on $\mathbb{R}^{n-k}$.

Denote by $T_{\Gamma}$ an $(n-k)$-dimensional torus
\begin{equation}\label{torus}
T_{\Gamma} = (e^{i\pi \langle \gamma_1,\varphi \rangle},...,e^{i\pi \langle \gamma_n,\varphi \rangle}) \subset \mathbb{C}^n,
\end{equation}
where $\varphi=(\varphi_1,...,\varphi_{n-k})\in \mathbb{R}^{n-k}$ and $\langle \gamma_j,\varphi \rangle=\gamma_{j,1}\varphi_1+...+\gamma_{j, n-k}\varphi_{n-k}$.

Consider a map
\begin{equation*}
\begin{gathered}
\widetilde{\psi}: \mathcal{R}\times T_{\Gamma} \rightarrow \mathbb{C}^n,\\
\widetilde{\psi}(u_1,...,u_n,\varphi) = (u_1e^{i\pi \langle \gamma_1,\varphi \rangle},...,u_ne^{i\pi \langle \gamma_n,\varphi \rangle}).
\end{gathered}
\end{equation*}
Let us show that $\widetilde{\psi}$ is not a nice map.  Let us define a group
\begin{equation*}
D_{\Gamma} = \Lambda^{*}/2\Lambda^{*}\approx \mathbb{Z}_2^{n-k}
\end{equation*}
and show that $D_{\Gamma}$ acts on $\mathcal{R}$. Let $\varepsilon \in D_{\Gamma}$ be a nontrivial element. We see that if $(u_1,...,u_n) \in \mathcal{R}$, then
\begin{equation*}
\varepsilon \cdot (u_1,...,u_n) = (u_1\cos\pi \langle \varepsilon,\gamma_1 \rangle,..,u_n\cos\pi \langle \varepsilon,\gamma_n \rangle)\in \mathcal{R}
\end{equation*}
because $\cos\pi \langle \varepsilon,\gamma_i \rangle = \pm 1$. We get that
\begin{equation*}
\widetilde{\psi}(u_1,...,u_n,\varphi) = \widetilde{\psi}(u_1\cos\pi \langle \varepsilon,\gamma_1 \rangle,..,u_n\cos\pi \langle \varepsilon,\gamma_n \rangle, \varphi+\varepsilon).
\end{equation*}
This shows that many points of $\mathcal{R} \times T_{\Gamma}$ have the same images under $\tilde{\psi}$. To fix this we need to take quotient of $\mathcal{R}\times{T_{\Gamma}}$ by the group $D_{\Gamma}$
\begin{equation}\label{mainmanifold}
\begin{gathered}
(u_1,...,u_n,\varphi) \sim (u_1\cos\pi \langle \varepsilon,\gamma_1 \rangle,..,u_n\cos\pi \langle \varepsilon,\gamma_n \rangle, \varphi+\varepsilon),\\
\mathcal{N} = \mathcal{R}\times_{D_{\Gamma}}T_{\Gamma}.
\end{gathered}
\end{equation}
The action of $D_{\Gamma}$ is free because it is free on the second factor. Therefore, $\mathcal{N}$ is a smooth $n$-manifold. So, we have a well-defined map
\begin{equation}\label{mainmap}
\begin{gathered}
\psi: \mathcal{N} \rightarrow \mathbb{C}^n,\\
\psi(u_1,...,u_n,\varphi) = (u_1e^{i\pi \langle \gamma_1,\varphi \rangle},...,u_ne^{i\pi \langle \gamma_n,\varphi \rangle}).
\end{gathered}
\end{equation}
As we noticed $D_{\Gamma}$ acts freely on $T_{\Gamma}$. Therefore, the projection
\begin{equation*}
\mathcal{N} = \mathcal{R}\times_{D_{\Gamma}} T_{\Gamma} \rightarrow T_{\Gamma}/D_{\Gamma} = T^{n-k}
\end{equation*}
onto the second factor defines a fiber bundle over $(n-k)$-torus $T_{\Gamma}/D_{\Gamma} = T^{n-k}$, where the fiber is $\mathcal{R}$.

Let us note that we need only system $(\ref{eqmain})$ to define $\psi$, $T_{\Gamma}$ and $\mathcal{N}$. We discussed in the previous section that there exists a polytope $P$ associated to system $(\ref{eqmain})$. We have the following theorem:

\begin{theorem}\label{mainemb} (see \cite{MironovCn}, \cite{MirPan}).
The Lagrangian $\psi(\mathcal{N})$ is immersed. The map $\psi$ is an embedding if and only if the polyhedron $P$ corresponding to system (\ref{eqmain}) is Delzant.
\end{theorem}

We see from ($\ref{mainmanifold}$) that the forms $d\varphi_1,...,d\varphi_{n-k}$ are closed invariant forms under the action of $D_{\Gamma}$, therefore they are elements of $H^1(\mathcal{N}, \mathbb{R})$. The lattice of the torus $T_{\Gamma}/D_{\Gamma} = T^{n-k}$ is formed by generators of $\Lambda^{*}$. Let $\varepsilon_1,..., \varepsilon_{n-k}$ be a basis for $\Lambda^{*}$. By $\varepsilon_{i, p}$ denote the $p$th coordinate of $\varepsilon_i$. Let $I_{\omega}$ be the symplectic area homomorphism and $I_{\mu}$ be the Maslov class of $\psi(\mathcal{N}) \subset \mathbb{C}^n$, where
\begin{equation*}
\omega = \frac{i}{2}\sum_{j=1}^{n}dz_j \wedge d\bar{z}_j.
\end{equation*}

Let us consider an $(n-k)$-dimensional vector
\begin{equation}\label{sympform}
\gamma_1+...+\gamma_n = (t_1,...,t_{n-k})^{T}.
\end{equation}

\begin{lemma}\label{hlem}( see \cite{Og} Section $5$).
The Maslov class $I_{\mu}$ of $\psi(\mathcal{N})$ is given by
\begin{equation}\label{lagrmasl}
I_{\mu} = t_1d\varphi_1+...+t_{n-k}d\varphi_{n-k}.
\end{equation}
If $\mathcal{R}$ is connected and simply connected, then
\begin{equation*}
H_1(\mathcal{N}, \mathbb{Z}) = H_1(T^{n-k}) = \mathbb{Z}^{n-k},
\end{equation*}
and there exists a basis $r_1,...,r_{n-k}$ for $H_1(\mathcal{N}, \mathbb{Z})$ such that
\begin{equation}\label{area}
\begin{gathered}
I_{\omega}(r_i) = \frac{\pi}{2} \sum\limits_{p=1}^{n-k}\varepsilon_{i, p}\delta_p, \\
I_{\mu}(r_i) = \sum\limits_{p=1}^{n-k}\varepsilon_{i,p}t_p
\end{gathered}
\end{equation}
where $\delta_p$ are defined in (\ref{eqmain}).

\end{lemma}

\textbf{Remark.} The formulas for the Maslov class and the symplectic area hold true when $\mathcal{R}$ is not simply connected. We do not want to go into details because in this paper we consider only Lagrangians with simply connected $\mathcal{R}$.
\\

Using the lemma above we can prove the following theorem:

\begin{theorem}\label{Fanoemb}(\cite{Og}).
\emph{Let $P$ be a Delzant and irredundant polytope and $L \subset \mathbb{C}^n$ be the corresponding Lagrangian submanifold. Then $L$ is monotone if and only if the polytope $P$ is Fano.}
\end{theorem}

\section{Proofs of Theorem $\ref{rest1}$ and Theorem $\ref{exist1}$}\label{proofs}

First, let us proof Theorem $\ref{rest1}$. We assumed that $n, p$ are even numbers and it's obvious that $2$ divides either $n$, or $p$. So, without loss of generality we can assume that $N_L \geqslant 3$. This allows us to work with lifted Floer homology.  The universal cover $\tilde{L}$ of $L$ is diffeomorphic to $S^{p-1} \times S^{n-p-1}$. As it is discussed in Section $\ref{prelim}$ we have the spectral sequence which converges to $FL^{\tilde(L)}(L) = 0$ with $E_1^{t,s} = H_{t+s - tN_L}(\tilde{L}, \mathbb{Z}_2) \otimes A^{p}$, $d_1 = \delta_1 \otimes T^{-1}$, where
\begin{equation*}
\delta_1: H_{m}(\tilde{L}, \mathbb{Z}_2) \rightarrow H_{m - 1 + N_L }(\tilde{L}, \mathbb{Z}_2).
\end{equation*}
We have that $E^{0,0}_1 = H_0(S^{p-1} \times S^{n-p-1}, \mathbb{Z}_2) = \mathbb{Z}_2$. We know that the spectral sequence converges to $0$, and therefore $\delta_r(H_0)$ should be nonzero for some $r$. there are only three possible cases:
\\
\\
1. $\delta_r(H_0(\tilde{L}, \mathbb{Z}_2)) = H_{p-1}(\tilde{L}, \mathbb{Z}_2))=\mathbb{Z}_2$
\\
\\
2. $\delta_r(H_0(\tilde{L}, \mathbb{Z}_2)) = H_{n-p-1}(\tilde{L}, \mathbb{Z}_2))=\mathbb{Z}_2$
\\
\\
3. $\delta_r(H_0(\tilde{L}, \mathbb{Z}_2)) = H_{n-2}(\tilde{L}, \mathbb{Z}_2))=\mathbb{Z}_2$
\\
\\
In the first case we have $0-1 +rN_L = p-1$ and $rN_L = p$. This means that $N_L$ divides $p$. In the second case we get that $rN_L = n-p$ and $N_L$ divides $n-p$. In the third case we have $rN_L = n-1$. The third case never happens because the left hand side is even (L is orientable) and the right hand side is odd (we assumed that $n$ is even). So, we proved that $N_L$ divides either $p$, or $n-p$. Theorem $\ref{rest1}$ is proved.
\\

Let us prove Theorem $\ref{exist1}$. Monotone Lagrangian embeddings of $S^{p-1} \times S^{n-p-1} \times T^2$ with minimal Maslov number $gcd(p, n-p+k)$ were constructed by the first author in \cite{Og}. Let us explain main ideas of the construction. Assume that $P $ is defined by inequalities
\begin{equation*}
\begin{gathered}
\left\{
 \begin{array}{l}
x_i + 1 \geqslant 0 \quad i=1,...,p-1 \\
-x_1 - ... - x_{p-1} +1 \geqslant 0 \\
x_i + 1 \geq 0 \quad i=p,...,n-2\\
- x_{p} - ... - x_{n-2} + 1 \geqslant 0  \\
 \end{array}
\right.
\end{gathered}
\end{equation*}
Easy to see that the system of inequalities above defines $P  = \Delta^{p-1} \times \Delta^{n-p-1}$, where $\Delta^{p-1}$, $\Delta^{n-p-1}$ are simplicies of dimensions $p-1$ and $n-p-1$, respectively. The corresponding system of quadrics has the form
\begin{equation}\label{eq}
\begin{gathered}
\left\{
 \begin{array}{l}
u_1^2 +...+u_p^2 = p\\
u_{p+1}^2 + ... + u_n^2 = n-p  \\
 \end{array}
\right.
\end{gathered}
\end{equation}
and defines $S^{p-1} \times S^{n-p-1}$.
\\

Let us add a parameter $k$. Let $P_k$ be a polytope defined by
\begin{equation*}
\begin{gathered}
\left\{
 \begin{array}{l}
x_i + 1 \geqslant 0 \quad i=1,...,p-1 \\
-x_1 - ... - x_{p-1} +1 \geqslant 0 \\
x_i + 1 \geq 0 \quad i=p,...,n-2\\
-x_1 - ... - x_{k} - x_{p} - ... - x_{n-2} + 1 \geqslant 0  \\
 \end{array}
\right.
n-p+k > p, \; \; k < p-1
\end{gathered}
\end{equation*}
Denote by $L_k$ the Lagrangian submanifold of $\mathbb{C}^n$ associated to the polytope $P_k$.  We see that the polytope $P_k$ is Delzant and Fano. Hence, from Theorem $\ref{mainemb}$ and Theorem $\ref{Fanoemb}$ we get that $L_k$ is embedded monotone Lagrangian.  The corresponding system of quadrics has the following form:
\begin{equation}\label{eq}
\begin{gathered}
\left\{
 \begin{array}{l}
u_1^2 +...+u_p^2 = p\\
u_1^2 + ... + u_{k}^2 + u_{p+1}^2 + ... + u_n^2 = n-p + k \\
 \end{array}
\right.
n - p +k > p, \;\; k < p-1
\end{gathered}
\end{equation}
The system above is equivalent to
\begin{equation}\label{dendr}
\begin{gathered}
\left\{
 \begin{array}{l}
2u_1^2 + ... + 2u_k^2 + u_{k+1}^2 + ... + u_n^2 = n+k\\
(n-2p+k)u_1^2 + ... + (n-2p+k)u_k^2 + (n-p+k)u_{k+1}^2 +  ... + \\
(n-p+k)u_p^2 -  pu_{p+1}^2 - ... - pu_n^2 = 0
 \end{array}
\right.  \\
\end{gathered}
\end{equation}
The second equation of the system defines a cone over the product of two ellipsoids of dimensions $p - 1$ and $n-p-1$. By intersecting it with ellipsoid of dimension $n-1$, defined by the first equation, we obtain that the system $(\ref{eq})$ defines
\begin{equation*}
\mathcal{R} = S^{p-1} \times S^{n-p-1} \subset \mathbb{R}^n.
\end{equation*}

From $(\ref{eq})$ we have that the corresponding torus $T_{\Gamma}$ and the embedding of $\mathcal{R} \times_{D_{\Gamma}} T_{\Gamma}$ into $\mathbb{C}^n$ are given by
\begin{equation*}
\begin{gathered}
\gamma_1 = ... = \gamma_k = (1,1)^T, \quad \gamma_{k+1}=...=\gamma_p = (1,0)^T, \quad \gamma_{p+1}=...=\gamma_n = (0,1)^T, \\
T_{\Gamma} = (\underbrace{e^{i\pi(\varphi_1 + \varphi_2)},...,e^{i\pi(\varphi_1 + \varphi_2)}}_k, \underbrace{e^{i\pi\varphi_1},...,e^{i\pi\varphi_1}}_{p-k},\underbrace{e^{i\pi\varphi_2},...,e^{i\pi\varphi_2}}_{n-p}) \subset \mathbb{C}^n,
\end{gathered}
\end{equation*}
\begin{equation}\label{emb1}
\begin{gathered}
\psi(u_1,...,u_n, \varphi_1, \varphi_2) = \\
 (u_1e^{i\pi(\varphi_1 + \varphi_2)}, ...,u_ke^{i\pi(\varphi_1 + \varphi_2)}, u_{k+1}e^{i\pi\varphi_1},...,u_pe^{i\pi\varphi_1},u_{p+1}e^{i\pi\varphi_2},...,u_ne^{i\pi\varphi_2})\\
\varphi_1, \varphi_2 \in \mathbb{R}.
\end{gathered}
\end{equation}

The dual lattice $\Lambda^{*}$ is generated by $\varepsilon_2 = (1,0)$, $\varepsilon_3 = (0,1)$. We know that $\mathcal{R} \times_{D_{\Gamma}} T_{\Gamma} \rightarrow T^2 = T_{\Gamma}/D_{\Gamma}$ is a fibration, where the fiber is $\mathcal{R}$. We have the fibration over $T^2$ with fiber $S^{p-1} \times S^{n-p-1}$ and with transition maps $\varepsilon_2, \varepsilon_3 \in D_{\Gamma}$
\begin{equation}\label{trans}
\begin{gathered}
\varepsilon_2(u_1,...,u_n ) \rightarrow  (-u_1,...,-u_{p}, u_{p+1},...,u_n ),\\
\varepsilon_3(u_1,...,u_n) \rightarrow  (-u_1,...,-u_{k}, u_{k+1},...,u_{p}, -u_{p+1},...,-u_n).
\end{gathered}
\end{equation}
If $k,p,n$ are even numbers, then
\begin{equation*}
\begin{gathered}
(-u_1\cos\phi - u_2\sin\phi, u_1\sin\phi - u_2\cos\phi,...,\\
-u_{p-1}\cos\phi - u_{p}\sin\phi, u_{p-1}\sin\phi - u_{p}\cos\phi, u_{p+1},..., u_n),\\
 \phi \in [0,\pi]
\end{gathered}
\end{equation*}
belongs to $\mathcal{R}$ for all $\phi$. The expression above defines the isotopy between $\varepsilon_2|_{fiber}$ and identity map. Similarly,  $\varepsilon_3|_{fiber}$ is isotopic to the identity map. Therefore, our fiber bundle is trivial.

In general the fibration is orientable if and only if numbers $p$ and $n-p+k$ are even (see \cite{Og}, Section 6).

From Formula $(\ref{area})$ and $(\ref{eq})$ we obtain
\begin{equation*}
I_{\mu}(r_1) = p, \quad I_{\mu}(r_2) = n-p+k, \quad I_{\omega}(r_1) = \frac{\pi}{2}p, \quad I_{\omega}(r_2) = \frac{\pi}{2}(n-p +k).
\end{equation*}
Hence, $L_k$ is monotone Lagrangian with minimal Maslov number $gcd(p,n-p+k)$ and monotonicity constant $\frac{\pi}{2}$.
\\
\\

Let us prove the second part of Theorem $\ref{exist1}$. We want to prove that if $n \geqslant 2p$, then varying $k$ we can get any divisor of $p$ as $\gcd\{p, n-p+k\}$. Let $D_{p}$ be the set of positive even divisors of $p$. We need to prove the following:
$$
\lbrace\gcd(p, n-p+k)\mid k\in 2\mathbb{Z}_{\geq 0}, k\leq p-2  \rbrace = D_p.
$$
We need to show that the right hand side is contained in the left hand side. Assume that $d\in D_{p}, dl=p, \; d\in 2\mathbb{Z}_{\geq 0}$. Then there exists a unique integer $m$ such that $d(ml+1)\geqslant n-p$ and $d(ml-l+1)<n-p$. If we take $k=d(ml+1)-(n-p)$, then $k \leqslant p-2$ and $\gcd(p, n-p+k)=\gcd(dl, d(ml+1))=d$.

\section{Proof of Theorem $\ref{rest2}$}

Let us prove the first part of the theorem. We assumed that $p$ is even and this means that $2$ divides $p$. So, without loss of generality we can assume that $N_L \geqslant 3$. The universal cover $\tilde{L}$ of $L$ is diffeomorphic to $(S^{p-1})^m$. Our Lagrangian $L$ belongs to $\mathbb{C}^{m(p-1)+l}$, hence $FH^{\tilde{L}}(L) = 0$. For simplicity we  use notation $H_{*}$ for $H_{*}((S^{p-1})^m, \mathbb{Z}_2)$. We have the following inequalities:
\begin{equation*}
\begin{gathered}
dim(V^{t,s}_2) = dim (\frac{Ker(\delta_1: V^{t,s}_1 \rightarrow V^{t-1,s}_1)}{Im(\delta_1: V^{t+1,s}_1 \rightarrow V^{t,s}_1)}) \geqslant dim(V^{t,s}_1) - dim(V^{t-1,s}_1) - dim(V^{t+1,s}_1) =
\\
dim(H_{t+s - tN_L}) - dim(H_{t+s - 1 - (t-1)N_L}) - dim(H_{t+s +1 - (t+1)N_L}).
\end{gathered}
\end{equation*}

Using the inequality above and the fact $dim(V^{t,s}_2) \leqslant dim(V^{t,s}_1)$ we see that
\begin{equation*}
\begin{gathered}
dim(V^{t,s}_3) \geqslant dim(V^{t,s}_2) - dim(V^{t-2,s+1}_2) - dim(V^{t+2,s-1}_2) \geqslant \\
dim(V^{t,s}_1) - dim(V^{t-1,s}_1) - dim(V^{t+1,s}_1) - dim(V^{t-2,s+1}_1) - dim(V^{t+2,s-1}_1),
\end{gathered}
\end{equation*}

\vspace{.16in}

\begin{equation*}
\begin{gathered}
0= dim(V^{t,s}_g) \geqslant dim(V^{t,s}_{g-1}) - dim(V^{t-g,s+g - 1}_{g-1}) - dim(V^{t+g,s-g+1}_{g-1}) \geqslant \\
dim(V^{t,s}_{1}) - \sum\limits_{i=1}^g dim(V^{t-i,s+i-1}_{1}) - \sum\limits_{i=1}^g dim(V^{t+i,s-i+1}_{1}) = \\
dim(H_{t+s - tN_L}) - \sum\limits_{j=1}^g dim(H_{t+s - 1 - (t-j)N_L}) - \sum\limits_{j=1}^g dim(H_{t+s + 1 - (t+j)N_L})
\end{gathered}
\end{equation*}

Suppose that $t=0$ and $s = [\frac{m}{2}](p-1)$. We have

\begin{equation*}
\sum\limits_{j=1}^g dim(H_{[\frac{m}{2}](p-1) - 1 + jN_L}) + \sum\limits_{j=1}^g dim(H_{[\frac{m}{2}](p-1) + 1 - jN_L}) \geqslant dim(H_{ [\frac{m}{2}](p-1) }).
\end{equation*}

\textbf{Assume that $N_L$ does not divide $p$.} Note that if $ - 1 + rN_L = p-1$, then $rN_L = p$ and $N_L$ divides $p$.  Also, $- 1 + rN_L \neq 2k(p-1)$ because the left hand side is odd and the right hand side is even. So, from the inequality above we have
\begin{equation*}
dim(H_{[\frac{m}{2}](p-1)}) \leqslant  \sum\limits_{j=1}^{[\frac{g}{2}]+1} dim(H_{[\frac{m}{2}](p-1) - 1 + 2jN_L}) + \sum\limits_{j=1}^{[\frac{g}{2}]+1} dim(H_{[\frac{m}{2}](p-1) + 1 - 2jN_L})
\end{equation*}
We see that $dim(H_{k}) = \binom{m}{r}$ if $k = r(p-1)$ and $dim(H_{k}) = 0$ otherwise. We get
\begin{equation*}
\begin{gathered}
\binom{m}{[\frac{m}{2}]} = dim(H_{ [\frac{m}{2}](p-1) }) \leqslant
\\
\sum\limits_{j=1}^{[\frac{g}{2}]+1} dim(H_{[\frac{m}{2}](p-1) - 1 + 2jN_L}) + \sum\limits_{j=1}^{[\frac{g}{2}]+1} dim(H_{[\frac{m}{2}](p-1) + 1 - 2jN_L}) \leqslant
\\
\binom{m}{0} + \binom{m}{1} + \binom{m}{3} +...+\binom{m}{[\frac{m}{2}]-3} +
\\
 \binom{m}{[\frac{m}{2}]+3} + \binom{m}{[\frac{m}{2}]+5} +...+ \binom{m}{m} <
\\
\sum\limits_{i=0}^{[\frac{m}{2}] - 3}\binom{m}{i} + \sum\limits_{i=[\frac{m}{2}] + 3}^{m}\binom{m}{i}.
\end{gathered}
\end{equation*}

As a result we have
\begin{equation}\label{expr}
\binom{m}{[\frac{m}{2}]} = dim(H_{ [\frac{m}{2}](p-1) }) \leqslant \sum\limits_{i=0}^{[\frac{m}{2}] - 3}\binom{m}{i} + \sum\limits_{i=[\frac{m}{2}] + 3}^{m}\binom{m}{i}.
\end{equation}

\begin{lemma}
$\binom{m}{[\frac{m}{2}]} > \sum\limits_{i=0}^{[\frac{m}{2}] - 3}\binom{m}{i} + \sum\limits_{i=[\frac{m}{2}] + 3}^{m}\binom{m}{i}  $.
\end{lemma}

Let us assume for a moment that the lemma is proved. Note that the lemma contradicts to $(\ref{expr})$. The inequality in the lemma is true because there are no terms $\binom{m}{[\frac{m}{2}] - 1}$ and $\binom{m}{[\frac{m}{2}] + 1}$. So, we have $ - 1 + jN_L = p-1$ for some $j$. Therefore, $jN_L =p $ and we get that $N_L$ divides $p$.
\\
\\
\emph{Proof of the lemma.} Let us assume for simplicity that $m$ is even. Almost the same proof works when $m$ is odd. We have
\begin{equation*}
2^m = (1+1)^m = \binom{m}{0} + ... + \binom{m}{m}.
\end{equation*}
Denote by
\begin{equation*}
\begin{gathered}
R_1 = \binom{m}{0} +...+\binom{m}{\frac{m}{2} - 3} + \binom{m}{\frac{m}{2} + 3} + ... + \binom{m}{m},\\
R_2 = \binom{m}{\frac{m}{2} - 2} + \binom{m}{\frac{m}{2} - 1} + \binom{m}{\frac{m}{2} } + \binom{m}{\frac{m}{2} +1} + \binom{m}{\frac{m}{2} +2}, \\
2^m = R_1 + R_2.
\end{gathered}
\end{equation*}
We need to show that
\begin{equation*}
\begin{gathered}
R_1 < \binom{m}{\frac{m}{2} } = R_2 - \binom{m}{\frac{m}{2} - 2} - \binom{m}{\frac{m}{2} - 1} - \binom{m}{\frac{m}{2} + 1} - \binom{m}{\frac{m}{2} +2} = \\
R_2 - 2\binom{m}{\frac{m}{2} + 1} - 2\binom{m}{\frac{m}{2} + 2} \\
 \Leftrightarrow \\
2^m - R_2 < R_2 - 2\binom{m}{\frac{m}{2} + 1} - 2\binom{m}{\frac{m}{2} + 2} \\
\Leftrightarrow \\
2^{m-1} < R_2 - \binom{m}{\frac{m}{2} + 1} - \binom{m}{\frac{m}{2} + 2} = \binom{m}{\frac{m}{2}} + \binom{m}{\frac{m}{2}-1} + \binom{m}{\frac{m}{2}-2} =\\
 \binom{m}{\frac{m}{2}} + \binom{m}{\frac{m}{2}+1} + \binom{m}{\frac{m}{2}+2}
\end{gathered}
\end{equation*}
So, we need to prove that $2^{m-1} < \binom{m}{\frac{m}{2}} + \binom{m}{\frac{m}{2}+1} + \binom{m}{\frac{m}{2}+2}$. Let us prove by induction the stronger inequality
\begin{equation}\label{ineq}
2^{m-1} < \binom{m}{\frac{m}{2}} + \binom{m}{\frac{m}{2}+1}.
\end{equation}
When $m=4$ inequality $\ref{ineq}$ is true. Suppose $\ref{ineq}$ is true for $m$ and let us prove that the inequality holds true for $m+1$. We have
\begin{equation*}
2^{m-1} < \binom{m}{\frac{m}{2}} + \binom{m}{\frac{m}{2}+1} = \binom{m+1}{\frac{m}{2}} = \frac{1}{2}\binom{m+1}{\frac{m}{2}} + \frac{1}{2}\binom{m+1}{\frac{m}{2} + 1}.
\end{equation*}
Since $[\frac{m+1}{2}] = \frac{m}{2}$ we have $2^m < \binom{m+1}{\frac{m}{2}} + \binom{m+1}{\frac{m}{2} + 1}$.

\section{Proof of Theorem $\ref{exist2}$}\label{proofexist2}

Let $P_k \in \mathbb{R}^{n-2}$ be a polytope defined by inequalities
\begin{equation*}
\begin{gathered}
\left\{
 \begin{array}{l}
x_i + 1 \geq 0 \quad i=1,...,n-2 \\
-x_1 - ... - x_{n-2} +1 \geq 0 \\
-x_1-...-x_k + k+2 \geq 0\\
 \end{array}
\right.  \\
0\leq k < n-1, \; \;\; k > \frac{n-3}{2}, \;\;\; n>3, \;\;\; k \;\; is \;\; even, \;\;\; n \;\; is \;\; odd
\end{gathered}
\end{equation*}
Let us note that under our conditions the polytope $P_k$ is redundant, i.e. the last inequality can be removed without changing the polytope. Indeed, let us sum inequalities $x_j + 1 \geqslant 0$ and $-x_1 -..-x_{n-2} + 1 \geqslant 0$, where $j=k+1,...,n-2$. We get $-x_1 -...-x_k + n-k-1 \geqslant 0$. We asumed that $2k > n-3$ and this implies that $n-k-1 < k+2$. Therefore, $-x_1-..-x_k \geqslant -(n-k-1) > -(k+2)$. As a result we get that under our conditions the last inequality is always satisfied.

Let $L_k$ be the Lagrangian submanifold of $\mathbb{C}^n$ associated to polytope $P_k$. If we remove the last inequality of $P_k$, then we get $(n-2)-$simplex defined by
\begin{equation*}
\begin{gathered}
\left\{
 \begin{array}{l}
x_i + 1 \geq 0 \quad i=1,...,n-2 \\
-x_1 - ... - x_{n-2} +1 \geq 0 \\
 \end{array}
\right.
\end{gathered}
\end{equation*}
Let us denote the simplex above by $\tilde{P}$. Let $\mathcal{R}_{P_k}$ and $\mathcal{R}_{\tilde{P}}$  be the submanifolds associated to polytopes $P_k$ and $\tilde{P}$, respectively. From Theorem $\ref{redundant}$ we see that $\mathcal{R}_{P_k}$ is diffeomorphic to $\mathcal{R}_{\tilde{P}} \times \mathbb{Z}_2$, i.e. $\mathcal{R}_{P_k}$ is disjoint union of two copies  of $\mathcal{R}_{\tilde{P}}$. The system of quadric associated to $\tilde{P}$ has the form
\begin{equation*}
u_1^2 + ... + u_{n-1}^2 = n-1
\end{equation*}
and associated to $P_k$
\begin{equation}\label{dissystem}
\begin{gathered}
\left\{
 \begin{array}{l}
u_1^2 + ... + u_{n-1}^2 = n-1\\
u_1^2 + ... + u_k^2 +  u_n^2 = 2k+2
 \end{array}
\right.
\end{gathered}
\end{equation}
We see that $\mathcal{R}_{\tilde{P}}$  is diffeomorphic to $S^{n-2}$ and $\mathcal{R}_{P_k}$ is diffeomorphic to $S^{n-2} \times \mathbb{Z}_2$. From $(\ref{dissystem})$ we have
\begin{equation*}
\gamma_1=...=\gamma_k = (1,1)^T, \quad \gamma_{k+1}=...=\gamma_{n-1} = (1,0)^T, \quad \gamma_{n} = (0,1)^T.
\end{equation*}
The lattice $\Lambda$ is generated by vectors $\gamma_{n-1}$, $\gamma_n$. The dual lattice $\Lambda^{*}$ is generated by $\varepsilon_{n-1} = (1,0)$ and $\varepsilon_n = (0, 1)$. The group $D_{\Gamma} = \Lambda^{*}/2\Lambda^{*}$ is isomorphic to $\mathbb{Z}_2^2$. Then the corresponding torus $T_{\Gamma}$ is given by
\begin{equation*}
T_{\Gamma} = (\underbrace{e^{i\pi(\varphi_1 + \varphi_2)},...,e^{i\pi(\varphi_1 + \varphi_2)}}_k, \underbrace{e^{i\pi\varphi_1},...,e^{i\pi\varphi_1}}_{n-k-1},e^{i\pi\varphi_2}) \subset \mathbb{C}^n.
\end{equation*}

From $(\ref{mainmap})$ we have embedding of $\mathcal{R}_{P_k} \times_{D_{\Gamma}} T_{\Gamma}$ into $\mathbb{C}^n$ given by
\begin{equation}\label{emb1}
\begin{gathered}
\psi(u_1,...,u_n, \varphi_1, \varphi_2) = \\
 (u_1e^{i\pi(\varphi_1 + \varphi_2)}, ...,u_ke^{i\pi(\varphi_1 + \varphi_2)}, u_{k+1}e^{i\pi\varphi_1},...,u_{n-1}e^{i\pi\varphi_1},u_{n}e^{i\pi\varphi_2})\\
\varphi_1, \varphi_2 \in \mathbb{R}.
\end{gathered}
\end{equation}
The Lagrangian associated to $P_k$ is $L_k = \psi(\mathcal{R}_{P_k} \times_{D_{\Gamma}} T_{\Gamma})$. First, let us note that $P_k$ is Delzant and by Theorem $\ref{mainemb}$ the Lagrangian $L_k$ is embedded.

We know that $\mathcal{R} \times_{D_{\Gamma}} T_{\Gamma}$ fibers over $T_{\Gamma}/D_{\Gamma} = T^2$ with fiber $S^{n-2} \times S^{0}$ and with transition maps $\varepsilon_{n-1}, \varepsilon_n \in D_{\Gamma}$
\begin{equation*}
\begin{gathered}
\varepsilon_{n-1}(u_1,..,u_n ) \rightarrow  (-u_1,...,-u_{n-1},u_n, ),\\
\varepsilon_n(u_1,...,u_n, ) \rightarrow (-u_1,...,-u_{k}, u_{k+1},...,u_{n-1},-u_n,).
\end{gathered}
\end{equation*}
We assumed that $n$ is odd and $k$ is even. Hence, we have isotopy between $\varepsilon_{n-1}|_{fiber}$ and the identity map defined by
\begin{equation*}
\begin{gathered}
(-u_1\cos\phi - u_2\sin\phi, u_1\sin\phi - u_2\cos\phi,...,\\
-u_{n-2}\cos\phi - u_{n-1}\sin\phi, u_{n-2}\sin\phi - u_{n-1}\cos\phi, \;\; u_n),\\
 \phi \in [0,\pi]
\end{gathered}
\end{equation*}

Also,
\begin{equation*}
\begin{gathered}
(-u_1\cos\phi - u_2\sin\phi, u_1\sin\phi - u_2\cos\phi,...,\\
-u_{k-1}\cos\phi - u_{k}\sin\phi, u_{k-1}\sin\phi - u_{k}\cos\phi, \;\; u_{k+1},...,u_{n-1}, -u_n),\\
 \phi \in [0,\pi]
\end{gathered}
\end{equation*}
defines isotopy between $\varepsilon_n|_{fiber}$  and
\begin{equation}\label{isoglue}
(u_1,...,u_n) \rightarrow (u_1,...,u_{n-1},-u_n).
\end{equation}
So, $\varepsilon_{n-1}$ is isotopic to identity map and we obtain that
\begin{equation*}
L_k = (S^{n-2} \times S^0 \times S^1) \times_{\varepsilon_n} S^1.
\end{equation*}
We can think about $L_k$ as $(S^{n-2} \times S^0 \times S^1) \times [0,1]$  with glued boundaries $(S^{n-2} \times S^0 \times S^1) \times \{0\}$ and $(S^{n-2} \times S^0 \times S^1) \times \{1\}$ by map $\varepsilon_n$. Easy to see from $(\ref{isoglue})$ that $\varepsilon_n$ joins copies of $S^{n-2} \times S^{0}$. So, we get that $L_k$ is diffeomorphic to $S^{n-2}\times T^2$.

The fiber is not connected and we can not use Lemma $\ref{hlem}$. Let us find the Maslov class and symplectic area form directly.

We assumed that $n>3$, therefore $\pi_1(\mathcal{R}) = 0$. As we discussed, $L_k$ fibers over $2-$torus with fiber $S^{n-2} \times S^0$ ( we proved that the fibration is trivial). The exact sequence of the fibration has the form
\begin{equation}\label{esequence}
0 \rightarrow \pi_1(L_k) \rightarrow \pi_1(T^2) \rightarrow \pi_0(\mathcal{R}) \rightarrow 0,
\end{equation}
where $T^2 = T_{\Gamma}/D_{\Gamma}$. Note that $\pi_0(\mathcal{R})$ is nontrivial, hence $\pi_1(L_k) \rightarrow \pi_1(T^2)$ is not surjective.
Let us consider a path $e_1: [0,1] \rightarrow L_k$
\begin{equation*}
e_1(s) =  \psi(u_1\cos(\pi s),...,u_{n-1}\cos(\pi s), u_n, s, 0), \;\;\; u_1,...,u_n = const,
\end{equation*}
and a path $e_2: [0,1] \rightarrow L_k$
\begin{equation*}
e_2(s) =  \psi(u_1,...,u_n, 0, 2s), \;\;\; u_1,...,u_n = const.
\end{equation*}
Note that
\begin{equation*}
\begin{gathered}
e_1(1) = \psi(-u_1,...,-u_{n-1}, u_n, 1, 0) =  \psi(\varepsilon_{n-1}(u_1,...,u_n, 0, 0)) = \psi(u_1,...,u_n, 0, 0) = e_1(0)
\end{gathered}
\end{equation*}
and this means that $e_1(s)$ is a loop. From $(\ref{mainmap})$ we get that $e_2(s)$ is a loop too.

Denote by $\sigma$ the projection $L_k \rightarrow T^2$. The torus $T_{\Gamma}/ D_{\Gamma} = T^2 \subset \mathbb{R}^2$ is formed by vectors $\varepsilon_{n-1}, \varepsilon_{n}$ and $\pi_1(T_{\Gamma}/ D_{\Gamma})$ is generated by $r_1 = t\varepsilon_{n-1}$ and $r_2 = t\varepsilon_n$, where $t \in [0,1]$. Then we see that $\sigma_{*}(e_1) = r_1$ and $\sigma_{*}(e_2) = 2r_2$ as elements of $\pi_1(T_{\Gamma}/ D_{\Gamma})$. Let $j$ be the map $\pi_1(T^2) \rightarrow \pi_0(\mathcal{R})$ from the exact sequence $(\ref{esequence})$. We have that loops $e_1$, $e_2$ generate $\pi_1(L_k)$ because $j(r_2) \neq 0$.

So, we proved that $H_1(L_k, \mathbb{Z})$ is generated by the cycles $e_1$, $e_2$. Let us find the Maslov class $I_{\mu_k}$ of $L_k$. From $(\ref{lagrmasl})$ we have
\begin{equation*}
\begin{gathered}
I_{\mu_k} = (n-1)d\varphi_1 + (k+1)d\varphi_2,\\
I_{\mu_k}(e_1) = n-1, \quad I_{\mu_k}(e_2)  = 2k+2.
\end{gathered}
\end{equation*}
 Let us also define a cycle
\begin{equation*}
\widetilde{e_1}(s) =  \psi(u_1,...,u_n, 2s, 0), \;\;\; u_1,...,u_n = const.
\end{equation*}
We see that $\sigma_{*}(\widetilde{e_1}) = 2r_1$  and $\widetilde{e_1} = 2e_1$ in $H_1(L_k, \mathbb{Z})$.

We consider $\mathbb{C}^n$ with the standard symplectic structure
\begin{equation*}
\omega = \sum\limits_{j=1}^{n}dx_j \wedge dy_j, \quad z_j = x_j + iy_j
\end{equation*}
and the Liouville form is given by
\begin{equation*}
\lambda = \sum\limits_{j=1}^n x_jdy_j.
\end{equation*}
From $(\ref{mainmap})$ we have
\begin{equation*}
x_j = u_j\cos(\pi<\gamma_j,\varphi>), \quad y_j = u_j\sin(\pi<\gamma_j,\varphi>),
\end{equation*}
\begin{equation*}
\begin{gathered}
\psi^{*}(x_jdy_j) = \frac{1}{2}\sin(2\pi<\gamma_j,\varphi>)u_jdu_j + \pi\sum\limits_{p=1}^{n-k} u_j^2\gamma_{j,p}\cos^2(\pi<\gamma_j,\varphi>)d\varphi_p =
\\
\frac{1}{4}\sin(2\pi<\gamma_j,\varphi>)d(u_j^2) + \pi\sum\limits_{p=1}^{n-k} u_j^2\gamma_{j,p}\cos^2(\pi<\gamma_j,\varphi>)d\varphi_p.
\\
\psi^{*}(\lambda) = \frac{1}{4}\sum\limits_{j=1}^{n}\sin(2\pi<\gamma_j,\varphi>)d(u_j^2) + \pi\sum\limits_{p=1}^{n-k} \delta_p\cos^2(\pi<\gamma_j,\varphi>)d\varphi_p.
\end{gathered}
\end{equation*}
As a result we obtain
\begin{equation*}
\begin{gathered}
I_{\omega}(e_1) = \frac{1}{2}I_{\omega}(\widetilde{e_1}) = \frac{1}{2}\lambda(\widetilde{e_1}) = \frac{\pi(n-1)}{2}, \\
I_{\omega}(e_2) = \lambda(e_2) = \pi(k+1).
\end{gathered}
\end{equation*}

We see that $L_k$ is monotone Lagrangian with minimal Maslov number
\begin{equation*}
N_{L_k} = gcd\{n-1, 2k+2\}.
\end{equation*}

\vspace*{.11in}

Let us prove that ``half of all'', or ``all but one'' divisors of $n-1$ can be realized. We assumed that $\frac{n-3}{2} <k \leqslant n-2$, which is equivalent to $n-1 < 2k+2 \leqslant  n-2$. Let $d$ be an even positive divisor of $n-1$.
\\

1. Assume $n-1$ is divisible by $4$. If $d \equiv 2$ mod $4$, then $d=4l+2$ for some non negative integer $l$. We have $n-1 < n-1+d \leqslant 2n - 2$. So, we see that $n-1+d \equiv 2$ mod 4 and $n-1+d = 2k+2$ for some even $k$ and
\begin{equation*}
4l+2 = d = gcd(n-1, n-1+d) = gcd(n-1, 2k+2), \quad \frac{n-3}{2} <k \leqslant n-2.
\end{equation*}

2. Assume $n-1 \equiv 2$ mod $4$. If $d \equiv 0$ mod $4$, then $n-1 < n-1+d \leqslant 2n-2$. We have that $n-1+d = 2k+2$ for some $k$ and
\begin{equation*}
d = gcd(n-1, n-1+d)= gcd(n-1, 2k+2), \quad \frac{n-3}{2} <k \leqslant n-2.
\end{equation*}
If $d \equiv 2$ mod $4$, then $n-1 = dl$ for some $l \geqslant 3$. Then $n-1 < n-1+2d < n-1 +dl = 2n-2$. We have that $n-1+2d = 2k + 2$ for some $k$ and
\begin{equation*}
d = gcd(n-1, n-1+d) = gcd(n-1, 2k+2), \quad \frac{n-3}{2} <k \leqslant n-2.
\end{equation*}

\section{Proof of Theorem $\ref{exist3}$}

The first part of the theorem is proved in \cite{Og} (Theorem 1.5). Let us prove the second part. The universal cover $\tilde{L}$ of $L$ is diffeomorphic to $\#_5 (S^{2p-1} \times S^{3p-2})$.  Homologies of $\tilde{L}$ are
\begin{equation*}
\begin{gathered}
H_0 = \mathbb{Z}_2, H_{2p-1} = \mathbb{Z}_2^5, H_{3p-2} = \mathbb{Z}_2^5, H_{5p-3} = \mathbb{Z}_2.
\end{gathered}
\end{equation*}
Arguing as in the proof of Theorem $\ref{rest1}$ we see that $FH^{\tilde{L}}(L) = 0$ and the corresponding spectral sequence converges to $0$. The group $H_{3p-2}$ can be killed only by $H_{2p-1}$. So, we have $2p-1 -1 +rN_L = 3p-2$. We get $rN_L = p$ and this implies that $N_L$ divides $p$.

\bibliographystyle{amsplain}

\vspace*{.3in}

Vardan Oganesyan, \\
L.D. Landau Institute for Theoretical Physics, Chernogolovka, 142432, Russia \\
\emph{E-mail}: vardan8oganesyan@gmail.com
\\

Yuhan Sun, \\
Department of Mathematics, Stony Brook University, Stony Brook, New York,11794, USA \\
\emph{E-mail}:  sun.yuhan@stonybrook.edu

\end{document}